\documentclass[10pt,draft,reqno]{amsart}
     \makeatletter
     \def\section{\@startsection{section}{1}%
     \z@{.7\linespacing\@plus\linespacing}{.5\linespacing}%
     {\bfseries
     \centering }}
     \def\@secnumfont{\bfseries}
     \makeatother
\theoremstyle{definition}

\theoremstyle{remark}

\numberwithin{equation}{section}
\usepackage{eucal,hyperref}
\RequirePackage{bbold}

\begin{document}

\title[Markov Branching Process without finite second moment]
    {On conditioned limit structure of the Markov Branching Process without finite second moment}

\author{{Azam A.~Imomov}}
\address{Azam Abdurakhimovich Imomov 
\newline\hphantom{iii} Karshi State University,
\newline\hphantom{iii} 17, Kuchabag street,
\newline\hphantom{iii} 180100 Karshi city, Uzbekistan.}
\email{{{imomov{\_}\,azam@mail.ru}}}
\thanks{\copyright \ 2021 Imomov A.A}



\dedicatory{Dedicated to my son Imron}

\vspace{.4cm}

\begin{abstract}
    Consider the continuous-time Markov Branching Process.
    In critical case we consider a situation when the generating
    function of intensity of transformation of particles has the infinite
    second moment, but its tail regularly varies in sense of Karamata.
    First we discuss limit properties of transition functions of the
    process. We prove  local limit theorems and investigate ergodic
    properties of the process.  Further we investigate limiting
    probability function conditioned to be  never extinct. Hereupon
    we obtain a new stochastic population process as a continuous-time
    Markov chain called the Markov Q-Process. We study main
    properties of Markov Q-Process.
\end{abstract}

\maketitle

\section{Introduction and preliminaries}             \label{MySec:1}

    Introducing the population of monotype individuals that are capable to perish and
    transforms into individuals of random number of the same type, we are interested
    in its evolution. These individuals may be biological kinds, molecules in chemical
    reactions etc. The most primitive mathematical model of population growth initiated
    by famous English statisticians H.Watson and F.Galton (1874) which is called now the
    Galton-Watson process. Yule (1924), considering the birth-and-death process, studied
    an evolution of biologic individuals. Feller (1939) used this model in problem of
    "struggle for existence". The Feller's problem was discussed by Neyman (1956), (1961)
    in situation of epidemic spread. The birth-and-death process was also studied by
    D'Ancona (1954) and Kendall (1948a), (1948b). In the book of Bharucha-Reid (1960) applications
    of models of particles evolution processes with Markov properties in the physics and
    biology were discussed. Kolmogorov and Dmitriev (1947) considered a population process
    which is an extension on the continuous-time case of definition of the Galton-Watson
    process and called the Markov Branching Process (MBP).

    Letting $Z(t)$ be the population size at the moment $t \in {\mathcal T}= [0;\, + \infty )$
    in MBP, we have the homogeneous continuous-time Markov chain $\left\{ {Z(t), t
    \in {\mathcal T}} \right\}$  with the state space ${\mathcal S}_0 = \{ 0\}  \cup {\mathcal S}$,
    where ${\mathcal S} \subset  \mathbb{N}=\{ 1,2, \ldots \} $. Evolution of the process
    occurs by the following scheme. Each individual existing at epoch $t \in {\mathcal T}$,
    independently of his history and of each other for a small time interval
    $(t;\,t + \varepsilon )$ transforms into $j \in {\mathcal S}_0 \backslash \{ 1\} $
    individuals with probability $a_j \varepsilon  + {o}(\varepsilon )$ and,
    with probability $1 + a_1 \varepsilon  + {o}(\varepsilon )$ each individual
    survives or makes evenly one descendant (as $\varepsilon  \downarrow 0$).
    Here $\left\{ {a_j } \right\}$ are intensities of individuals' transformation that
    $a_j  \ge 0$ for $j \in {\mathcal S}_0 \backslash \{ 1\} $ and
$$
    0 < a_0  <  - a_1  = \sum\limits_{j \in {\mathcal S}_0 \backslash \{ 1\} } {a_j }  < \infty.
$$
    Appeared new individuals undergo transformations under same way as above.
    The Markovian nature of the process yields that its transition functions
$$
    P_{ij} (t) = \mathbb{P}_i \left\{ {Z(t) = j} \right\}:
    = \mathbb{P}\left\{ {Z(t + \tau ) = j\left| {Z(\tau ) = i} \right.} \right\}
$$
    satisfy the Kolmogorov-Chapman equation
$$
    P_{ij} (t) = \sum\limits_{k \in {\mathcal S}_0 } {P_{ik} (u)
    \cdot P_{kj} (t - u)}
    \quad \parbox{2cm}{\textit{for} {} $u \le t$,}   \eqno (1.1)
$$
    and the following branching property holds for all $i,j \in {\mathcal S}_0 $:
$$
    P_{ij} (t) = \sum\limits_{j_1  + j_2  +  \cdots  + j_i  = j}
    {P_{1j_1 } (t) \cdot P_{1j_2 } (t)\, \cdots \, P_{1j_i } (t)};   \eqno (1.2)
$$
    see Athreya and Ney (1972, Ch. III). Thus, for studying of evolution of
    $\left\{ {Z(t), t \in {\mathcal T}} \right\}$ is suffice to set the transition
    functions $P_{1j} (t)$. These functions in turn, can be expressed using the
    local densities $\left\{ {a_j } \right\}$ by relation
$$
    P_{1j} (\varepsilon ) = \delta _{1j}  + a_j \varepsilon  + o(\varepsilon )
    \quad \parbox{2cm}{\textit{as} {} $\varepsilon  \downarrow 0$,}       \eqno (1.3)
$$
    where $\delta _{ij} $ represents Kronecker's delta function. A probability
    generating  function (GF) is a main analytical tool in our discussions on MBP.
    A GF version of relation (1.3) is
$$
    F(\tau ;s) = s + f(s) \cdot \tau  + o(\tau )
    \quad \parbox{2cm}{\textit{as} {} $\tau  \downarrow 0,$}
$$
    for all $0 \le s < 1$, where $F(t;s) = \sum\nolimits_{j \in {\mathcal S}_0}{P_{1j}(t)s^j}$
    and $f(s) = \sum\nolimits_{j \in {\mathcal S}_0 } {a_j s^j }$.

    GF $F(t;s)$ satisfies to the functional equation
$$
    F(t + \tau ;s) = F\left( {t;F(\tau ;s)} \right),   \eqno (1.4)
$$
    for any $t, \tau  \in {\mathcal T}$ with the boundary condition $F(0;s) = s$.
    Moreover it satisfies the equation
$$
    {{\partial F(t;s)} \over {\partial t}} = f\left( {F(t;s)} \right),   \eqno (1.5)
$$
    the backward Kolmogorov equation. It follows from theory of
    differential equations that the solution of this equation
    represents unique GF which satisfies the
    equation (1.4); see Athreya and Ney (1972, p.106).

    If the offspring mean $a: = \sum\nolimits_{j \in {\mathcal S}} {ja_j }= f'(s \uparrow 1)$
    is finite then $F(t;1) = 1$; see Asmussen and Hering (1983, p.119).
    By using equations (1.1), (1.2) and (1.5) it can be computed that
    $\mathbb{E}_i Z(t) := \sum\nolimits_{j \in {\mathcal S}_0 }{jP_{ij} (t)}  = ie^{at} $.
    The last formula shows that long-term properties of MBP seem variously depending
    on value of parameter $a$. Hence the MBP is classified as critical if $a = 0$
    and sub-critical or supercritical if $a < 0$ or $a > 0$ respectively.

    Further we write everywhere $\mathbb{P}\{ * \} $ and
    $\mathbb{E}[*]$ instead of $\mathbb{P}_1 \{*\} $ and
    $\mathbb{E}_1 [*]$ respectively.

    Let a random variable ${\mathcal H} = \inf \left\{ {t \in {\mathcal T}:Z(t) = 0} \right\}$
    be the extinction time of MBP. The fundamental extinction theorem states that
    $\mathbb{P}_i \left\{ {{\mathcal H}<\infty} \right\}=q^i$, where
    $q = \inf \left\{ {s \in (0,1]:f(s) = 0} \right\}$ is the extinction
    probability that $q=1$ if the process is non-supercritical. An asymptote of
    probability of ${\mathcal H}$ has first been observed by Sevastyanov (1951).
    Exertions of this variable were treated also by Heatcote \textit{et al.} (1967),
    Nagaev and Badalbaev (1967), Zolotarev (1957). Put the conditioned distribution
$$
    \mathbb{P}_i^{{\mathcal H}(t)} \{  * \} : = \mathbb{P}_i
    \left\{ { * \left| {t < {\mathcal H} < \infty } \right.} \right\}.
$$
    Sevastyanov (1951) proved that in the sub-critical case there exits
    a limiting distribution law
    $P_j^ *   = \lim _{t \to \infty } \mathbb{P}^{{\mathcal H}(t)}
    \left\{ {Z(t) = j} \right\}$ with GF
$$
    \sum\limits_{j \in {\mathcal S}} {P_j^ *  s^j }  = 1 - \exp
    \left\{ {a\int_0^s {{{dx} \over {f(x)}}} } \right\},  \eqno (1.6)
$$
    if and only if $\sum\nolimits_{j \in {\mathcal S}} {a_j j\ln j}  < \infty $.
    In the critical situation he also proved that if $2b: = f''(s \uparrow 1)$
    is finite, then ${{Z(t)} \mathord{\left/ {\vphantom {{Z(t)} {bt}}} \right.
    \kern-\nulldelimiterspace} {bt}}$ has a limiting exponential law.

    In the discrete-time situation $\mathbb{P}_i^{{\mathcal H}(t + \tau )}\{*\}$
    converge as $\tau  \to \infty $ to a probability measure, which defines
    homogeneous Markov chain called the Q-process; see Athreya and Ney (1972, pp.56--60).
    The Q-process was considered first by Lamperti and Ney (1968). Some properties
    of it were discussed by Pakes (1971, 1999, 2010),
    Imomov (2001, 2002, 2014b, 2014c).

    Similarly in the MBP case a limit
    $\lim _{\tau  \to \infty } \mathbb{P}_i^{{\mathcal H}(t + \tau )}
    \left\{ {Z(t) = j} \right\}$ has an honest probability measures
    $\left\{ {{\mathcal Q}_{ij} (t)} \right\}$ which defines
    a homogeneous continuous-time process as a Markov chain
    and, called the Markov Q-Process; see Imomov (2012).
    Let $W(t)$ be the state size at the moment $t \in {\mathcal T}$ in Markov
    Q-Process. Then $W(0)\mathop  = \limits^d Z(0)$ and
    $\mathbb{P}_i \left\{ {W(t) = j} \right\} = {\mathcal Q}_{ij} (t)$. In the
    paper Imomov (2012) some asymptotic properties of the chain
    $\left\{ {W(t),t \in {\mathcal T}} \right\}$ are observed. Namely it was
    proved that if the associated MBP is critical and $f''(s \uparrow 1)$
    is finite, then ${{W(t)} \mathord{\left/ {\vphantom {{W(t)}
    {\mathbb{E}W(t)}}} \right. \kern-\nulldelimiterspace} {\mathbb{E}W(t)}}$ has
    a limiting Erlang's law. In this case there is an invariant measure.
    In the non-critical situation under the condition when (1.6) holds,
    an invariant distribution exists for the process
    $\left\{ {W(t),t \in {\mathcal T}} \right\}$.

    In this paper we consider MBP without further power moments.
    In non-critical case we rest satisfied only with the condition
    $\sum\nolimits_{j \in {\mathcal S}} {a_j j\ln j}  < \infty $.
    In the critical case our reasoning is bound up with elements
    of the theory of regularly varying functions in sense of
    Karamata; see Karamata (1933). We remember that real-valued,
    positive and  measurable function $\ell(x)$ is  said to be
    slowly varying (SV) at $\alpha $ if it belongs to a class
$$
    \mathfrak{S}_\alpha := \left\{ {\ell(x) \in \mathbb{R}_+: \;
    \mathop {\lim }\limits_{x \to \alpha} {{\ell(\lambda x)} \over {\ell(x)}}
    = 1, \quad   \forall \lambda  \in \mathbb{R}_+ } \right\}.
$$
    A function $\textsf{V}(x)$ is said to be regularly varying (RV) at
    $\alpha $ with index of regular variation $\rho \in \mathbb{R}_+$
    if it in the form $\textsf{V}(x) = x^\rho \ell(x)$, where
    $\ell(x) \in \mathfrak{S}_\alpha $. We denote $\mathfrak{R}_\alpha^\rho$
    be the class of RV functions. It is evidently that
    $\mathfrak{S}_\alpha   \equiv \mathfrak{R}_\alpha ^0 $.

    Throughout the paper in the critical case we suppose that the
    infinitesimal GF $f(s)$ satisfies the condition
$$
    f(s) = (1 - s)^{1 + \nu } \mathfrak{L}
    \left( {{1 \over {1 - s}}} \right),   \eqno (1.7)
$$
    for $0 \le s < 1$, where $0 < \nu \le 1$ and
    $\mathfrak{L}(x) \in \mathfrak{S}_\infty$.
    By the assumed criticality
$$
    {{f''(s \uparrow 1)} \over 2} = \mathop {\lim }\limits_{s \uparrow 1} {{f(s)}
    \over {(1 - s)^2 }} = \mathop {\lim }\limits_{x \downarrow 0}
    {1 \over {x^{1 - \nu } }}\mathfrak{L}\left( {{1 \over x}} \right) = \infty.
$$
    If $f''(s\uparrow 1) < \infty $, then condition (1.7) holds with $\nu  = 1$ and
    $\mathfrak{L}(t) \to {{f''(s\uparrow1)} \mathord{\left/
    {\vphantom {{f''(1)} 2}} \right. \kern-\nulldelimiterspace} 2}$
    as $t \to \infty $. Thus our process contains a process
    with the finite second moment.

    Section 2 is devoted to auxiliary lemmas these will be useful for our
    purpose. First we take assertions about asymptotical decay of the functions
    $ F(t;s)$ and ${{\partial F(t;s)} \mathord{\left/ {\vphantom
    {{\partial F(t;s)} {\partial s}}} \right. \kern-\nulldelimiterspace}
    {\partial s}}$. One of important facts of this section is the
    Monotone ratio Lemma on limit properties of the transition
    functions $P_{ij} (t)$. In consequence of this lemma we take
    complete accounts on asymptotic behaviors of  states of MBP.
    In Section 3 we observe invariance properties of states of MBP.
    In non-critical case we  hold on to results
    from the paper of Imomov (20104a).

    In Section 4 we consider the Markov Q-Process and discuss properties
    concerning its construction and asymptotic properties of transition
    functions $\left\{ {{\mathcal Q}_{ij} (t)} \right\}$. In particular we
    compute the q-matrix and the GF version of the Kolmogorov backward
    equation implied by $\left\{ {{\mathcal Q}_{ij} (t)} \right\}$.
    We also observe ergodic properties of the Markov Q-Process.

\section{Auxiliary results}

    Let
$$
    R(t;s) := q - F(t;s).
$$

\textbf{Lemma~1.}
    \textit{The following assertions are true for all $0 \le s < 1$}.
\begin{itemize}
    \item    {\it Let $a \ne 0$. Then}
$$
    R(t;s) = {\mathcal A}(t;s) \cdot \beta ^t,   \eqno (2.1)
$$
    {\it where $\beta : = \exp \left\{ {f'(q)} \right\}< 1$, and }
$$
    {\mathcal A}(t;s) = (q - s)\exp \left\{ {\int_s^{F(t;s)} {\left[
    {{1 \over {u - q}} - {{f'(q)} \over {f(u)}}} \right]du} } \right\}.
$$
    \item   {\it Let $a = 0$. If the condition (1.7) is assumed, then}
$$
    R(t;s) \sim {{{\mathcal N}(t)} \over {(\nu t)^{{1 \mathord{\left/
    {\vphantom {1 \nu }} \right. \kern-\nulldelimiterspace} \nu }} }}
    \cdot \left( {1 + {{\mathcal{M}(s)} \over t}} \right)^{ - {1 \mathord{\left/
    {\vphantom {1 \nu }} \right. \kern-\nulldelimiterspace} \nu }},      \eqno (2.2)
$$
    {\it where }
$$
    {\mathcal N}(n) \cdot \mathcal{L}^{{1 \mathord{\left/ {\vphantom {1 \nu }}
    \right. \kern-\nulldelimiterspace}\nu}}\left({{{\bigl(\nu n \bigr)^{{1
    \mathord{\left/ {\vphantom {1 \nu }} \right. \kern-\nulldelimiterspace}
    \nu }} } \over  {{\mathcal N}(n)}}} \right)\buildrel {} \over
    \longrightarrow 1
    \quad \parbox{2.2cm}{\textit{as} {} $n  \rightarrow \infty$}    \eqno (2.3)
$$
    {\it and }
$$
    \mathcal{M}(s) = \int_1^{{1 \mathord{\left/ {\vphantom {1 {(1 - s)}}}
    \right. \kern-\nulldelimiterspace}{(1 - s)}}}
    {{{dx} \over {x^{1 - \nu }\mathfrak{L}(x)}}}.       \eqno (2.4)
$$
\end{itemize}

\textbf{\textit{Proof}.}
    Let's consider first the non-critical case. Multiplying by
    $f'(q) \cdot \left( {F(t;s) - q} \right)$ the both sides of (1.5) yields
$$
    {{dF(t;s)} \over {F(t;s) - q}} \cdot \left[ {1 -
    {{f\left( {F(t;s)} \right) - f'(q) \cdot \left( {F(t;s) - q} \right)}
    \over {f\left( {F(t;s)} \right)}}} \right] = f'(q) \cdot dt.
$$
    After integration on $[0;\,t] \subset {\mathcal T}$ it follows from this equality
$$
    {{R(t;s)} \over {R(0;s)}} = \beta ^t \exp \left\{ {\int_s^{F(t;s)}
    {\left[ {{1 \over {u - q}} - {{f'(q)} \over {f(u)}}} \right]du} } \right\},
$$
    hereinafter $\beta  = \exp \left\{ {f'(q)} \right\}$.
    Since $R(0;s)=q-s$, this relation we can write in the form of (2.1).

    To prove the second part we transform the backward Kolmogorov equation (1.5)
    to the following integral equation:
$$
    \int_s^{F(t;s)} {{{dx} \over {f(x)}}}  = t.   \eqno (2.5)
$$
    Further we rewrite equation (2.5) in form of
$$
    \int_0^{F(t;s)} {{{dx} \over {f(x)}}}  = t + \mathcal{M}(s),   \eqno (2.6)
$$
    where
$$
    \mathcal{M}(s) = \int_0^s {{{dx} \over {f(x)}}}.
$$
    Seeing the condition (1.7) and  denoting
    $1 - x = {1 \mathord{\left/ {\vphantom {1 u}} \right.
    \kern-\nulldelimiterspace} u}$ it follows from (2.6) that
$$
    \int\limits_1^{{1 \mathord{\left/ {\vphantom {1 {R(t;s)}}}
    \right. \kern-\nulldelimiterspace} {R(t;s)}}} {{{dx}
    \over {x^{1 - \nu } \mathfrak{L}(x)}}}  = t + \mathcal{M}(s),     \eqno (2.7)
$$
    herein the function $\mathcal{M}(s)$ becomes (2.4).

    Now we recall the following property of SV functions.
    If $\ell(x) \in \mathfrak{S}_\infty  $ remains locally
    bounded in $[A; + \infty )$ for some $A \in \mathbb{R}_+$, then
$$
    \int_A^x {u^\lambda  \ell(u)du}
    \sim {1 \over {\lambda  + 1}}x^{\lambda  + 1} \ell(x),    \eqno (2.8)
$$
    as $x \to \infty $ for $\lambda  >  - 1$; see Bingham \textit{et al.} (1987, p.26).

    Since the upper bound of the integral in left-hand side of (2.7) grows to
    infinity as $t \to \infty $, it is possible to use the property (2.8) and
    we can get the following asymptotic formula:
$$
    \int\limits_1^{{1 \mathord{\left/ {\vphantom {1 {R(t;s)}}} \right.
    \kern-\nulldelimiterspace} {R(t;s)}}}{{{dx} \over {x^{1 - \nu }
    \mathfrak{L}(x)}}} = {{1+o(1)} \over {\nu R^\nu  (t;s) \cdot
    \mathfrak{L}\left( {{1 \mathord{\left/ {\vphantom {1 {R(t;s)}}}
    \right. \kern-\nulldelimiterspace} {R(t;s)}}} \right)}}
    \quad \parbox{2.2cm}{\textit{as} {} $t  \to \infty$.}
$$
    Combining this formula with (2.7) yields
$$
    R(t;s) = {{\mathfrak{L}^{-{1 \mathord{\left/ {\vphantom {1 \nu }}
    \right. \kern-\nulldelimiterspace} \nu }} \left( {{1 \mathord{\left/
    {\vphantom {1 {R(t;s)}}} \right. \kern-\nulldelimiterspace} {R(t;s)}}}
    \right)} \over {(\nu t)^{{1 \mathord{\left/{\vphantom {1 \nu }} \right.
    \kern-\nulldelimiterspace} \nu }}}}\cdot {{1+o(1)}\over
    {\left({1 + \displaystyle{{\mathcal{M}(s)}\over t}}
    \right)^{{1 \mathord{\left/{\vphantom {1 \nu }} \right.
    \kern-\nulldelimiterspace} \nu }} }}
    \quad \parbox{2.2cm}{\textit{as} {} $t  \to \infty$.}     \eqno (2.9)
$$
    Letting $q(t): = R(t;0) = \mathbb{P}\left\{ {Z(t) > 0} \right\}$ and seeing
    that $\mathcal{M}(0) = 0$, we obtain from (2.9) that
$$
    q(t) = {{\mathfrak{L}^{ - {1 \mathord{\left/ {\vphantom {1 \nu }} \right.
    \kern-\nulldelimiterspace} \nu }} \left( {{1 \mathord{\left/
    {\vphantom {1 {q(t)}}} \right. \kern-\nulldelimiterspace}
    {q(t)}}} \right)} \over {(\nu t)^{{1 \mathord{\left/ {\vphantom
    {1 \nu }} \right. \kern-\nulldelimiterspace} \nu }} }}{\left(1+o(1)\right)}
    \quad \parbox{2.2cm}{\textit{as} {} $t  \to \infty$,}
$$
    the surviving probability of the process $\left\{ Z(t) \right\}$. On the
    other hand it is  clear that ${{R(t;s)} \mathord{\left/ {\vphantom {{R(t;s)}
    {q(t)}}} \right. \kern-\nulldelimiterspace} {q(t)}} \to 1$.  Hence introducing
    the function ${\mathcal N}(t)$ satisfying property (2.3), we can write formula (2.2).

    The Lemma is proved.

    Since $F(t;s) \to q$ as $t \to \infty $, we obtain
    ${\mathcal A}(t;s) \to {\mathcal A}(s)$, where
$$
    {\mathcal A}(s)=(q - s)\exp \left\{{\int_s^{q}{\left[{{1 \over {u - q}}
    - {{f'(q)} \over {f(u)}}} \right]du} } \right\}.   \eqno (2.10)
$$
    Therefore it follows from (2.1) that
$$
    \mathbb{P}\left\{ {n < {\mathcal H} < \infty } \right\}
    \sim {\mathcal A}(0) \cdot \beta ^t
    \quad \parbox{2.2cm}{\textit{as} {} $t  \to \infty$}
$$
    if and only if
$$
    \int_0^q {{{|\ln \beta | \cdot x - f(q - x)} \over {xf(q - x)}}dx}
    = \ln {q \over {{\mathcal A}(0)}} < \infty.     \eqno (2.11)
$$
    Integral in the left-hand side in (2.11) can be transformed to a form
$$
    q\cdot \int_0^1 {{{\hat {a} \cdot x - \widehat f(1 - x)}
    \over {x \widehat f(1 - x)}}dx},
$$
    where $\widehat f(s) = \sum\nolimits_{j \in {\mathcal S}_{0} } {\widehat a_j s^j }$,
    $\widehat a_j  = a_j q^{j - 1}$ and $\widehat a = \sum\nolimits_{j \in {\mathcal S}}
    {j\widehat a_j }  = f'(q) < 1$. It is known (see Sevastyanov (1951)) that
    convergence of last integral is equivalent to a convergence of the series
    $\sum\nolimits_{j \in {\mathcal S}} {\widehat a_j j\ln j}$.
    As $\widehat a_j  < a_j$, the condition
$$
    \sum\limits_{j \in {\mathcal S}} {a_j j\ln j}  < \infty  \eqno (2.12)
$$
    is sufficient to be satisfied the condition (2.11).

    Further consider the function ${{\partial R(t;s)} \mathord{\left/
    {\vphantom {{\partial R(t;s)} {\partial s}}} \right.
    \kern-\nulldelimiterspace} {\partial s}}$.

\textbf{Lemma~2.}
    \textit{The following assertions are true for all $0 \le s < 1$}.
\begin{itemize}
    \item    {\it If $a \ne 0$, then}
$$
    {{\partial R(t;s)} \over {\partial s}} = {{f'(q)} \over {f(s)}}
    {\mathcal A}(t;s)  \beta ^t \left( {1 + o(1)} \right)
    \quad \parbox{2.2cm}{\it \textit{as} {} $t \to \infty$,}     \eqno (2.13)
$$
    {\it where the function ${\mathcal A}(t;s)$ is defined in (2.1).}

    \item   {\it Let $a = 0$. If the condition
    (1.7) is satisfied, then}
$$
    {{\partial R(t;s)} \over {\partial s}} =  - \left( {{{R(t;s)}
    \over {1 - s}}} \right)^{1 + \nu } {{\mathfrak{L}\left(
    {{1 \mathord{\left/ {\vphantom {1 {R(t;s)}}} \right.
    \kern-\nulldelimiterspace} {R(t;s)}}} \right)} \over
    {\mathfrak{L}\left( {{1 \mathord{\left/ {\vphantom
    {1 {(1 - s)}}} \right. \kern-\nulldelimiterspace}
    {(1 - s)}}} \right)}}.           \eqno (2.14)
$$
\end{itemize}

\textbf{\textit{Proof}.}
    Transform the backward Kolmogorov equation (1.5) to (2.5) and
    differentiating it with respect $s$, we have
$$
    {{\partial R(t;s)} \over {\partial s}}
    = - {{f\left( {F(t;s)} \right)} \over {f(s)}}.   \eqno (2.15)
$$
    In non-critical case $f(s) \sim f'(q)(s - q)$ as $s \to q$.
    Therefore (2.15) entails
$$
    {{\partial R(t;s)} \over {\partial s}} =
    {{f'(q)} \over {f(s)}}R(t;s)\left( {1 + o(1)} \right)
    \quad \parbox{2.2cm}{\it \textit{as} {} $t \to \infty$,}
$$
    for all $0 \le s < 1$. Using therein (2.1) we will get (2.13).

    Same way, for critical case, from (2.15) we obtain the following relation:
$$
    {{\partial R(t;s)} \over {\partial s}} =  - {{R^{1 + \nu }
    (t;s)} \over {f(s)}}\mathfrak{L}\left( {{1 \over {R(t;s)}}} \right).
$$
    This relation together with condition (1.7) produces (2.14).

    The proof of Lemma 2 is completed.

    One can see that the Lemma 2 has a simple appearance,
    but as it will be visible further, this lemma is important
    in our discussions. Namely it will easily be computed that
$$
    \left. {{{\partial F(t;s)} \over {\partial s}}} \right|_{s = 0}  = P_{11} (t),
$$
    the probability of return of the process with initial
    state $Z(0) = 1$ to the one through the time $t$.
    As $f(0) = a_0  > 0$, putting $s = 0$ in (2.13)
    and (2.14), we obtain the following local limit theorems.

\textbf{Theorem~1.}
    \textit{Let $a \ne 0$. Then}
$$
    \beta ^{-t}\cdot P_{11}(t)\sim {{|\ln \beta |}\over {a_0}} \mathcal{A}(t;0)
    \quad \parbox{2.2cm}{\it \textit{as} {} $t \to \infty$.}
$$
    \textit{If the condition (2.12) is satisfied,
    then $\mathcal{A}(t;0)\to \mathcal{A}(0)<\infty$ as $t \to \infty$.}

\textbf{Theorem~2.}
    \textit{Let $a = 0$. If the condition (1.7) is satisfied, then}
$$
    (\nu t)^{1 + {1 \mathord{\left/ {\vphantom {1 \nu }} \right.
    \kern-\nulldelimiterspace} \nu }}\cdot P_{11} (t)
    \sim {{{\mathcal N} (t)} \over {a_0 }}
    \quad \parbox{2cm}{\textit{as} {} $t \to \infty $,}    \eqno (2.16)
$$
    \textit{where the function ${\mathcal N}(t)$ satisfies the property (2.3).}

    Further we will use the following Monotone ratio limit property.

    \textbf{Lemma~3~(Imomov (2014a)).} \textit{For all $j \in {\mathcal S}$}
$$
    {{P_{1j} (t)} \over {P_{11} (t)}} \uparrow \mu _j  < \infty
    \quad \parbox{2.2cm}{\it \textit{as} {} $t \to \infty$.}     \eqno (2.17)
$$

    Now observe ${{P_{ij} (t)} \mathord{\left/ {\vphantom {{P_{ij} (t)}
    {P_{11} (t)}}} \right. \kern-\nulldelimiterspace} {P_{11} (t)}}$
    as $t \to \infty$ for $i,j \in {\mathcal S}$. We see that
$$
    {\textbf M}_{i} (t;s): = \sum\limits_{j \in {\mathcal S}}
    {{{P_{ij} (t)} \over {P_{11} (t)}}s^j } \longrightarrow
    iq^{i - 1}  \cdot {\textbf M}(s)
    \quad \parbox{2.2cm}{\it \textit{as} {} $t \to \infty$,}    \eqno (2.18)
$$
    for all $0 \le s < 1$, where
$$
    {\textbf M}(s) = \sum\limits_{j \in {\mathcal S}} {\mu _j s^j }.
$$
    By virtue of relation (2.18) to studying the long-term
    behavior of $P_{ij} (t)$ is suffice to consider the function
    ${\textbf M}(t;s): = {\textbf M}_{1} (t;s)$.
    The transition function version of (2.18) is
$$
    {{P_{ij} (t)} \over {P_{11} (t)}} \longrightarrow  iq^{i - 1} \mu _j
    \quad \parbox{2.2cm}{\it \textit{as} {} $t \to \infty$.}     \eqno (2.19)
$$
    Using (2.19) and from Theorems 1 and 2 we get complete accounts
    about asymptotic behaviors of transition function $P_{ij} (t)$.

\textbf{Theorem~3.}
    \textit{Let $a \ne 0$. If the condition (2.12) is satisfied, then}
$$
    \beta ^{- t} \cdot P_{ij} (t) \sim iq^{i - 1} \mu _j{{|\ln \beta |}
    \over {a_0 }}  \mathcal{A}(t;0)
    \quad \parbox{2.2cm}{\it \textit{as} {} $t \to \infty$.}
$$
    \textit{If the condition (2.12) is satisfied,
    then $\mathcal{A}(t;0)\to \mathcal{A}(0)<\infty$ as $t \to \infty$.}

\textbf{Theorem~4.}
    \textit{Let $a = 0$. If the condition (1.7) is satisfied, then}
$$
    (\nu t)^{1 + {1 \mathord{\left/ {\vphantom {1 \nu }} \right.
    \kern-\nulldelimiterspace} \nu }}\cdot P_{ij} (t)
    \sim {{{i \mu _j}} \over {a_0 }}{{\mathcal N} (t)}
    \quad \parbox{2cm}{\textit{as} {} $t \to \infty $,}
$$
    \textit{where the function ${\mathcal N}(t)$ satisfies the property (2.3).}

\section{Invariant properties of transition functions $P_{ij} (t)$}

    Continuing researches of the asymptote of transition functions
    $P_{ij} (t)$ we deal with problems of ergodicity and existence of
    invariant measure. Ergodicity properties of arbitrary continuous-time
    Markov chain are described in the monograph of Anderson (1991). The
    invariant (or stationary) measure of chain $\left\{ {Z(t),
    t \in {\mathcal T}} \right\}$ is a set of non-negative numbers
    $\left\{ {\nu _j ,j \in {\mathcal S}_0 } \right\}$ satisfying to the equation
$$
    \nu _j  = \sum\limits_{k \in {\mathcal S}_0 } {\nu _k P_{kj} (t)},    \eqno (3.1)
$$
    for any $t \in {\mathcal T}$. Equation (3.1) means an invariant
    property of the measure $\left\{ {\nu _j } \right\}$ concerning to
    the transition functions $\left\{ {P_{ij} (t)} \right\}$.
    If $\sum\nolimits_{j \in {\mathcal S}_0 } {\nu _j }  < \infty $
    (or without loss of generality $\sum\nolimits_{j \in {\mathcal S}_0 }
    {\nu _j }  = 1$) then it is called an invariant distribution.

    Further we will discuss the role of the set
    $\left\{ {\mu _j, j \in {\mathcal S}} \right\}$ defined in (2.17)
    as invariant measures. The following theorem holds.

\textbf{Theorem~5~(Imomov (2014a)).}
    \textit{Non-negative numbers $\left\{ {\mu _j } \right\}$
    satisfy the invariant equation}
$$
    \beta ^t  \cdot \mu _j
    = \sum\limits_{k \in {\mathcal S}} {\mu _k P_{kj} (t)},   \eqno (3.2)
$$
    \textit{for $j \in {\mathcal S}$ and for all $t \in {\mathcal T}$. The function}
    ${\textbf M}(s) = \sum\nolimits_{j \in {\mathcal S}} {\mu _j s^j }$
    \textit{satisfies the functional equation}
$$
    {\textbf M}\left( {F(t;s)} \right) = \beta ^t  \cdot
    {\textbf M}(s) + {\textbf M}\left( {F(t;0)} \right),  \eqno (3.3)
$$
    \textit{which converges for $0 \le s < 1$. Equation (3.3) has a unique
    solution that is power series with non-negative coefficients for $0 \le s < q$.}

    In the following two theorems the explicit forms
    of the function ${\textbf M}(s)$ will be obtained.

\textbf{Theorem~6~(Imomov (2014a)).}
    \textit{Let $a \ne 0$ and the condition (2.12) is satisfied. Then}
$$
    {\textbf M}(s) = {{a_0 } \over {\left| {f'(q)} \right|}} \cdot
    \left[ {1 - {{{\mathcal A}(s)} \over {{\mathcal A}(0)}}} \right].     \eqno (3.4)
$$

\textbf{Theorem~7.}
    \textit{Let $a = 0$. If the condition
    (1.7) is satisfied, then}
$$
    {\textbf M}(s) = a_0 \mathcal{M}(s),     \eqno (3.5)
$$
    \textit{where $\mathcal{M}(s)$ is form of (2.4).}

\textbf{\textit{Proof}.}
    Recall $R(t;s) = 1 - F(t;s)$ and write
$$
    {\textbf M}(t;s) = {{F(t;s) - F(t;0)} \over {P_{11} (t)}}
    = \left( {1 - {{R(t;s)} \over {q(t)}}} \right)
    \cdot {{q(t)} \over {P_{11} (t)}},      \eqno (3.6)
$$
    where $q(t): = R(t;0)$. Using the second part of the
    Lemma 1 and after elementary transformations we find
$$
    1 - {{R(t;s)} \over {q(t)}} = {{\mathcal{M}(s)}
    \over {\nu t}}\left( {1 + o(1)} \right)
    \quad \parbox{2.2cm}{\it \textit{as} {} $t \to \infty$.}     \eqno (3.7)
$$
    According to (2.16), ${{q(t)} \mathord{\left/ {\vphantom {{R(t)}
    {P_{11} (t)}}} \right. \kern-\nulldelimiterspace} {P_{11} (t)}}
    \sim a_0 \nu t$ as $t \to \infty $. Then considering relations
    (3.6) and (3.7) we get to (3.5).

    The theorem is proved.

    Previous two theorems along with Lemma 3 allows to judge about asymptotic
    behavior of the sum $\sum\nolimits_{j \in {\mathcal S}} {\mu _j }$.
    According to Lemma 3 this sum converges for $a < 0$ and diverges if $a > 0$.
    For the case $a = 0$ formulas (2.4), (2.8) and (3.5) show that
$$
    {\textbf M}(s) = {{a_0 } \over {\nu (1 - s)^\nu  }} \cdot
    \mathfrak{L}_{\mu } \left( {{1 \over {1 - s}}} \right)
    \quad \parbox{2cm}{\it \textit{as} {} $s\uparrow 1 $,}    \eqno (3.8)
$$
    where $\mathfrak{L}(x) \cdot \mathfrak{L}_{\mu }(x)\to 1$ as $x \to \infty$.
    The ensuing theorem follows from equality (3.8) according
    to the Hardy-Littlewood Tauberian theorem.

\textbf{Theorem~8.}
    \textit{Let $a = 0$. If the condition (1.7) is satisfied, then}
$$
    \sum\limits_{j = 1}^n {\mu _j }  = {{a_0 }
    \over {\nu ^2 \Gamma (\nu )}}n^\nu  \mathfrak{L}_{\mu } (n),
$$
    \textit{where $\Gamma ( * )$ is Euler's Gamma function and
    $\mathfrak{L}(x) \cdot \mathfrak{L}_{\mu }(x)\to 1$ as $x \to \infty$.}

    Theorem 3 shows that in non-critical situation transition functions
    $P_{ii} (t)$ have an exponential decay behavior as $t \to \infty $.
    The limit
$$
    \lambda _{\mathcal S} = - \mathop {\lim }\limits_{t \to \infty }
    {{\ln P_{ii} (t)} \over t}
$$
    independent on $i \in {\mathcal S}$ and characterizes a decay rate of state space
    of chain $\left\{ {Z(t)} \right\}$. It is called the decay parameter of states
    of the chain. MBP classified as $\lambda _{\mathcal S}$-transient if
    $\int_0^{+ \infty }{e^{\lambda _{\mathcal S}t} P_{ii}(t)dt} < \infty$
    and $\lambda _{\mathcal S} $-recurrent otherwise. In this case invariant measure
    is called $\lambda _{\mathcal S} $-invariant. According to the general classification
    MBP is called $\lambda _{\mathcal S} $-positive if $\lim _{t \to \infty }
    e^{\lambda _{\mathcal S} t} P_{ii} (t) > 0$ and $\lambda _{\mathcal S} $-null if this
    is zero; see Li \textit{et al.} (2010). Theorems 3 and 6 imply the following theorem.

\textbf{Theorem~9.}
    \textit{Let $a \ne 0$ and the condition (2.12) is satisfied.
    Then $\lambda _{\mathcal S}  = \left| {\ln \beta } \right|$ and the Markov
    chain $\left\{ {Z(t)} \right\}$ is $\lambda _{\mathcal S} $-positive.  The set
    $\left\{ {\mu _j, j \in {\mathcal S}} \right\}$ determined by  GF (3.4) is
    unique (up to multiplicative constant) $\lambda _{\mathcal S} $-invariant measure.}

    In critical case the set $\left\{ {\mu _i } \right\}$ directly enters to
    a role of invariant measure for MBP. Indeed, in this case $\beta  = 1$ and
    as it has been proved in Theorem 5 that
$$
    \mu _j  = \sum\limits_{k \in {\mathcal S}} {\mu _k P_{kj} (t)}, \;\; j \in {\mathcal S},
$$
    for all $t \in {\mathcal T}$.

    As shown in Theorems 3 and 4 hit probability of MBP to any state through
    the long interval time depends on the initial state. In other words an
    ergodic property is not carried out. Thereby we will seek an ergodic chain
    associated to MBP. Recalling ${\mathcal H}$ be the extinction moment of MBP we write
$$
    \mathbb{P}_i \left\{ {t < {\mathcal H} < \infty , Z(t) = j} \right\}
    = \mathbb{P}\left\{ {t < {\mathcal H}
    < \infty \left| {Z(t) = j} \right.} \right\} \cdot P_{ij} (t).
$$
    Since the probability of extinction of $j$ particles is $q^j $
    then it follows that
$$
    \mathbb{P}_i \left\{ {t < {\mathcal H} < \infty , Z(t) = j}
    \right\} = P_{ij} (t) \cdot q^j.       \eqno (3.9)
$$
    We have also that
$$
    \mathbb{P}_i \left\{ {t < {\mathcal H} < \infty } \right\} = \sum\limits_{j \in {\mathcal S}}
    {\mathbb{P}_i \left\{ {Z(t) = j,t < {\mathcal H} < \infty } \right\}}
    = \sum\limits_{j \in {\mathcal S}} {P_{ij} (t)q^j }.
$$
    Using the formula (3.9) from last relation we obtain that
$$
    \mathbb{P}_i \left\{ {t < {\mathcal H} < \infty } \right\}
    = \sum\limits_{j \in {\mathcal S}} {P_{ij} (t)q^j }.       \eqno (3.10)
$$

    Now put into consideration the conditional transition function
$$
    \mathbb{P}_i^{{\mathcal H}(t)} \{  * \} : = \mathbb{P}_i \left
    \{ { * \left| {t < {\mathcal H} < \infty } \right.} \right\}.
$$
    Let $\widetilde P_{ij}(t)= \mathbb{P}_i^{{\mathcal H}(t)} \left\{{Z(t)=j} \right\}$
    be a transition matrix which defines a new stochastic process
    $\left\{ {\widetilde Z(t), t \in {\mathcal T}} \right\}$. It is easy
    to be convinced that $\left\{ {\widetilde Z(t), t \in {\mathcal T}} \right\}$
    represents a homogeneous Markov chain. Indeed probabilities
    $\widetilde P_{ij} (t)$ satisfy to the Kolmogorov-Chapman equation (1.1) and
    have the branching property (1.2). According to last theorems properties
    of trajectories of $\left\{ {\widetilde Z(t),t \in {\mathcal T}}\right\}$
    lose dependence on the initial state as $t \to \infty $. Consider an
    appropriate GF
$$
    {\mathcal V}_i (t;s) = \sum\limits_{j \in {\mathcal S}} {\widetilde P_{ij} (t)s^j }.
$$

\textbf{Theorem~10~(Imomov (2014a)).}
    \textit{Let $a \ne 0$ and the condition
    (2.12) is satisfied. Then limits}
$$
    \mathop {\lim }\limits_{t \to \infty } \widetilde P_{ij} (t) = \nu _j
$$
    \textit{exist for all $i,j \in {\mathcal S}$ and these are determined by GF}
$$
    {\mathcal V}(s) = {{{\textbf M}(qs)} \over {{\textbf M}(q)}},       \eqno (3.11)
$$
    \textit{where function ${\textbf M}(s)$ is defined in (3.4).}

\textbf{Remark.}
    \textit{Theorem 10 is generalization of Sevastyanov's result (1.6) in which
    corresponding result established for the sub-critical situation only.
    Indeed it is easy to see that the limit probability GF (1.6)
    is the proprietary case of (3.11). The set $\left\{ {\nu _j } \right\}$
    represents a probability distribution. In fact setting $s = 1$ in (3.11) and
    taking into account equality (3.4), it follows that ${\mathcal V}(1)
    = \sum\nolimits_{j \in {\mathcal S}} {\nu _j }  = 1$. Moreover if the condition
    (2.12) is satisfied, then}
$$
    \sum\limits_{j \in {\mathcal S}} {j\widetilde P_{ij} (t)}
    \longrightarrow {q \over {{\mathcal A}(0)}}
    \quad \parbox{2.2cm}{\it \textit{as} {} $t \to \infty$,}
$$
    \textit{and ${\mathcal V}'(s \uparrow 1)={q \mathord{\left/{\vphantom
    {q {{\mathcal A}(0)}}}\right.   \kern-\nulldelimiterspace} {{\mathcal A}(0)}}$.}

    Under the condition of Theorem 10 for the MBP $\left\{ {\widetilde Z(t),
    t \in {\mathcal T}} \right\}$ exists the unique (up to multiplicative constant)
    set of non-negative numbers $\left\{ {\nu _i } \right\}$ which not all are
    zero and we see without difficulty that the GF ${\mathcal V}(s) = \sum\nolimits_{j
    \in {\mathcal S}} {\nu _j s^j }$ satisfies the invariance equation
$$
    \beta ^t  \cdot {\mathcal V}(s) = {\mathcal V}\left( {{{F(t;qs)}
    \over q}} \right) - {\mathcal V}\left( {{{F(t;0)} \over q}} \right).
$$
    So $\left\{ {\nu _i } \right\}$ is invariant measure.

    It follows from (3.9) and (3.10) that if $a \ne 0$ and
    the condition (2.12) is satisfied, then
$$
    \widetilde P_{ij} (t)={P_{ij} (t) \over
    {\sum\nolimits_{k \in {\mathcal S}} {P_{ik} (t)q^{k - j}}}}.
$$

    In the critical situation $\mathbb{P}\left\{ {{\mathcal H} < \infty } \right\} = 1$.
    Since $1 - F^i (t;s) \sim iR(t;s)$, we write
$$
    {\mathcal V}_i (t;s)\, \sim 1 - {{R(t;s)} \over {q(t)}}
    \quad \parbox{2.2cm}{\it \textit{as} {} $t \to \infty$,}      \eqno (3.12)
$$
    where $q(t)= R(t;0)$.

\textbf{Theorem~11.}
    \textit{Let $a = 0$. If the condition
    (1.7) is satisfied, then}
$$
    \nu t \cdot {\mathcal V}_i (t;s) \longrightarrow \mathcal{M}(s)
    \quad \parbox{2.2cm}{\it \textit{as} {} $t \to \infty$,}       \eqno (3.13)
$$
    \textit{where GF $\mathcal{M}(s)= \sum\nolimits_{j \in {\mathcal S}}
    {\textbf{m}_j s^j }$ is form of (2.4). Moreover}
$$
    \sum\limits_{j = 1}^n {\textbf{\textit{m}} _j }  = {{1}
    \over {\nu ^2 \Gamma (\nu )}}n^\nu  \mathfrak{L}_{\mu }(n),     \eqno (3.14)
$$
    \textit{where $\Gamma ( * )$ is Euler's Gamma function and
    $\mathfrak{L}(x) \cdot \mathfrak{L}_{\mu }(x)\to 1$ as $x \to \infty$.}

\textbf{\textit{Proof}} of the convergence (3.13) directly comes out from
    (3.7) and (3.12). Relation (3.14) follows from Theorem 7 and Theorem 8.

\section{The Markov Q-process}

    In this section we will be interested in a limiting interpretation of the
    conditioned transition function $\mathbb{P}_i^{{\mathcal H}(t + \tau )}
    \left\{ {Z(t) = j} \right\}$ as $\tau  \to \infty $ for all $t \in {\mathcal T}$.
    As it was said in first Section this limit is an honest probability measure.
    This measure defines so-called Markov Q-Process (MQP) be the continuous-time
    Markov chain $\left\{ {W(t), t \in {\mathcal T}} \right\}$ with the state space
    ${\mathcal E} \subset \mathbb{N}$. The random function $W(t)$ is the state size
    at the  moment $t \in {\mathcal T}$ in MQP. The transition function
    ${\mathcal Q}_{ij} (t) = \mathbb{P}_i \left\{ {W(t) = j} \right\}$ is form of
$$
    {\mathcal Q}_{ij} (t) = \mathop {\lim }\limits_{\tau  \to \infty }
    \mathbb{P}_i^{{\mathcal H}(t + \tau )} \left\{ {Z(t) = j} \right\}
    = {{jq^{j - i} } \over {i\beta ^t }}P_{ij} (t),     \eqno (4.1)
$$
    for $i,j \in {\mathcal E}$, where $\beta  = \exp \{ f'(q)\} $; see Imomov (2012).
    It is easy to be convinced that $0 < \beta  \le 1$ decidedly. To wit $\beta  = 1$
    if $a = 0$ and $\beta  < 1$ otherwise. In our presupposition the MBP is honest.
    Since $F\left( {t;q} \right) = q$ and $\left. {{{F\left( {t;s} \right)}
    \mathord{\left/ {\vphantom {{F\left( {t;s} \right)} {\partial s}}} \right.
    \kern-\nulldelimiterspace} {\partial s}}} \right|_{s = q}  = \beta ^t $, it follows
    from (4.1) that $\sum\nolimits_{j \in {\mathcal E}} {{\mathcal Q}_{ij} (t)}  = 1$.

    Combining equalities (1.3) and (4.1) we obtain the following representation:
$$
    {\mathcal Q}_{1j} (\varepsilon ) = \delta _{1j}  + \lambda_j \varepsilon  + o(\varepsilon ),
    \quad \parbox{1.6cm}{\textit{as} {} $\varepsilon  \downarrow 0$,}    \eqno (4.2)
$$
    with probability densities
$$
    \lambda_0  = 0, \quad \lambda_1  = a_1  - \ln \beta <0,
    \quad \mbox{\textit{and}} \quad \lambda_j  = jq^{j - 1} a_j  \ge 0
    \quad \parbox{3cm}{\textit{for} {} $j \in {\mathcal E}\backslash \{ 1\}$,}
$$
    where $\left\{ {a_j } \right\}$ are evolution intensities of MBP $Z(t)$.
    It follows from (4.2) that
$$
    g(s): = \sum\limits_{j \in {\mathcal E}} {\lambda_j s^j }
    = s\left[ {f'(qs) - f'(q)} \right].     \eqno (4.3)
$$
    Needles to see that this GF is infinitesimal one because $g(1) = 0$.
    So the infinitesimal GF $g(s)$ completely defines the process
    $W(t)$, where $\left\{ {\lambda_j } \right\}$ are intensities of process
    evolution satisfying $\lambda_j  > 0$ for $j \in {\mathcal E}\backslash \{ 1\} $ and
$$
    0 <  - \lambda_1  = \sum\limits_{j \in {\mathcal E}\backslash \{ 1\} } {\lambda_j }  < \infty.
$$

\subsection{Construction, existence and uniqueness}

    Let's now discuss basic properties of transition matrix
    $\mathbb{Q}(t) = \left\{ {{\mathcal Q}_{ij} (t)} \right\}$.
    Herewith we follow methods and facts from monograph of Anderson (1991).
    First we prove the following theorem.

\textbf{Theorem~12.}
    \textit{Let $\left\{ {W(t), t \in {\mathcal T}} \right\}$ be the MQP
    given by infinitesimal GF $g(s)$. Then the transition matrix
    $\mathbb{Q}(t)$ is standard and honest. Its components
    ${\mathcal Q}_{ij} (t)$ are positive and uniformly continuous
    function with respect to $t \in {\mathcal T}$ for all $i,j \in {\mathcal E}$.}

\textbf{\textit{Proof}.}
    According to the branching property (1.2), we see
$$
    P_{ij} (\varepsilon ) = \delta _{ij}  + ia_{j - i + 1}
    \varepsilon  + o(\varepsilon )
    \quad \parbox{1.6cm}{\textit{as} {} $\varepsilon  \downarrow 0$.}
$$
    Hence seeing equality (4.1) it follows
\begin{equation}
    \left\{\begin{array}{l} {\mathcal Q}_{ii} (\varepsilon )
    = 1 + \left( {ia_1  - \ln \beta } \right)\varepsilon
    + o(\varepsilon )\; ,  \\ \nonumber
    \\
    {\mathcal Q}_{ij} (\varepsilon ) = jq^{j - i} a_{j - i + 1}
    \varepsilon  + o(\varepsilon ) \; \hfill,   \nonumber
    \end{array} \right. \quad \parbox{2cm}{\textit{as}
    {} $\varepsilon  \downarrow   0$,}                  \eqno (4.4)
\end{equation}
    for all $i,j \in {\mathcal E}$. Considering representations (4.4) we have
\begin{eqnarray}
    \sum\limits_{j \in {\mathcal E}} {\left| {{\mathcal Q}_{ij} (\varepsilon )
    - \delta _{ij} } \right|} \nonumber
    & = & \sum\limits_{j \in {\mathcal E}\backslash \{ i\} }  {{\mathcal Q}_{ij}
    (\varepsilon )}  + \left| {{\mathcal Q}_{ii} (\varepsilon ) - 1} \right| \\ \nonumber
    & = & \sum\limits_{j \in {\mathcal E}\backslash \{ i\} }
    {{\mathcal Q}_{ij} (\varepsilon )}  + 1 - {\mathcal Q}_{ii} (\varepsilon )
    \le 2\left| {1 - {\mathcal Q}_{ii} (\varepsilon )} \right| \longrightarrow 0,  \nonumber
\end{eqnarray}
    as $\varepsilon  \downarrow 0$. So that ${\mathcal Q}_{ij} (t)$ is standard.
    A positiveness of functions ${\mathcal Q}_{ij} (t)$ is obvious owing to (4.4).
    The Markovian nature of the process $\left\{ {W(t)} \right\}$ implies
    the Kolmogorov-Chapman equation:
$$
    {\mathcal Q}_{ij} (t + \varepsilon ) = \sum\limits_{k \in {\mathcal E}}
    {{\mathcal Q}_{ik} (t){\mathcal Q}_{kj} (\varepsilon )}.
$$
    Hence supposing $\varepsilon  > 0$ it follows that
\begin{eqnarray}
    {\mathcal Q}_{ij} (t + \varepsilon ) - {\mathcal Q}_{ij} (t) \nonumber
    & = & \sum\limits_{k \in {\mathcal E}} {{\mathcal Q}_{ik}
        (\varepsilon ){\mathcal Q}_{kj} (t)}  - {\mathcal Q}_{ij} (t) \\ \nonumber
    & = & \sum\limits_{k \in {\mathcal E}\backslash \{ i\} }
        {{\mathcal Q}_{ik} (\varepsilon ){\mathcal Q}_{kj} (t)}
        - {\mathcal Q}_{ij} (t) \cdot \left[ {1 - {\mathcal Q}_{ii}
        (\varepsilon )} \right].  \nonumber
\end{eqnarray}
    It follows from here that
\begin{eqnarray}
    - \left[ {1 - {\mathcal Q}_{ii} (\varepsilon )} \right] \nonumber
    & \le &  - {\mathcal Q}_{ij} (t) \cdot \left[ {1 - {\mathcal Q}_{ii}
        (\varepsilon )} \right] \le {\mathcal Q}_{ij} (t + \varepsilon )
        - {\mathcal Q}_{ij} (t) \\ \nonumber
    & \le &  \sum\limits_{k \in {\mathcal E}\backslash \{ i\} }
        {{\mathcal Q}_{ik} (t){\mathcal Q}_{kj} (\varepsilon )}
        \le \sum\limits_{k \in {\mathcal E}\backslash \{ i\} }
        {{\mathcal Q}_{kj} (\varepsilon )}
        = 1 - {\mathcal Q}_{ii} (\varepsilon ),  \nonumber
\end{eqnarray}
    so $\left| {{\mathcal Q}_{ij} (t + \varepsilon )- {\mathcal Q}_{ij}(t)}
    \right| \le 1 - {\mathcal Q}_{ii} (\varepsilon )$.  Similarly
\begin{eqnarray}
    \left| {{\mathcal Q}_{ij} (t-\varepsilon )-{\mathcal Q}_{ij}(t)} \right|  \nonumber
    & = & \left| {{\mathcal Q}_{ij} (t) - {\mathcal Q}_{ij}
        (t - \varepsilon )} \right| \\ \nonumber
    & \le &  1 - {\mathcal Q}_{ii} \left( {t - (t - \varepsilon )}
        \right) = 1 - {\mathcal Q}_{ii} (\varepsilon ).   \nonumber
\end{eqnarray}
    Therefore we obtain $\left| {{\mathcal Q}_{ij} (t + \varepsilon )
    - {\mathcal Q}_{ij} (t)} \right| \le 1 - {\mathcal Q}_{ii}
    \left( {|\varepsilon |} \right)$ for any $\varepsilon \ne 0$.
    This relation implies that ${\mathcal Q}_{ij} (t)$ is uniformly
    continuous function with respect to $t \in {\mathcal T}$  because
    $\lim _{\varepsilon  \downarrow 0} {\mathcal Q}_{ii} (\varepsilon ) = 1$.

    The theorem is proved.

    It can easily be verified that a GF version of (4.4) is
$$
    G_i (t;s): = \mathbb{E}_i s^{W(t)}  = \sum\limits_{j \in {\mathcal E}}
    {{\mathcal Q}_{ij} (t)s^j }  = {{qs} \over {i\beta ^t }}\left[ {{\partial
    \over {\partial x}}\left( {{{F(t;x)} \over q}} \right)^i } \right]_{x = qs},
$$
    or more obviously
$$
    G_i (t;s) = \left[ {{{F(t;qs)} \over q}} \right]^{i - 1} G(t;s),     \eqno (4.5)
$$
    where
$$
    G(t;s) = G_1 (t;s) = {s \over {\beta ^t }}\left.
    {{{\partial F(t;x)} \over {\partial x}}} \right|_{x = qs} .
$$

\textbf{Theorem~13.}
    \textit{All states of the Markov chain $\left\{ {W(t)} \right\}$ are stable.
    The transition functions $\left\{ {{\mathcal Q}_{ij} (t)} \right\}$ are the
    Feller functions. These functions are differentiable and has a finite and
    continuous derivative with respect to $t \in {\mathcal T}$. Its q-matrix
    $\left\{ {q_{ij}  = {\mathcal Q}'_{ij} (\varepsilon  \downarrow 0)} \right\}$ has components}
\begin{equation}
    q_{ij}  =\left\{\begin{array}{l} i\lambda_1  + (i - 1)\ln \beta \; \hfill,
    \qquad \parbox{2.4cm}{\textit{when} {} $i = j $,}  \\ \nonumber
    \\
    \displaystyle {{j\lambda_{j - i + 1}} \over {j - i + 1}} \hfill,
    \qquad  \parbox{2.4cm}{\textit{when} {} $i \ne j $,}   \nonumber
    \end{array} \right.     \eqno (4.6)
\end{equation}
    \textit{where $\lambda_i $ are in (4.2) and $q_{ij}  \ge 0$
    when $i \ne j$ and, $q_i : = q_{ii}  < 0$ for all $i,j \in {\mathcal E}$.
    Moreover it satisfies the identity}
$$
    {\mathcal Q}'_{ij} (t + \tau ) = \sum\limits_{k \in {\mathcal E}}
    {{\mathcal Q}'_{ik} (\tau ){\mathcal Q}_{kj} (t)},
    \quad \parbox{3.2cm}{\textit{for any} {} $t, \tau  \in {\mathcal T} $,}    \eqno (4.7)
$$
    \textit{the backward Kolmogorov system.}

\textbf{\textit{Proof}.}
    It follows from the relation (4.4) that for all $i \in {\mathcal E}$
$$
    q_i = \mathop {\lim }\limits_{\varepsilon  \downarrow 0} {{1 - {\mathcal Q}_{ii}
    (\varepsilon )} \over {- \varepsilon }} <  + \infty,
$$
    that is all states are stable and also the right-sided derivative
    ${\mathcal Q}'_{ij} (\varepsilon  \downarrow 0)$ is finite.

    From the relation (4.5) we have
$$
    {\mathcal Q}_{ij} (t) < {{G(t;s)} \over {s^j }} \cdot \left[ {\widehat {F}(t;s)} \right]^{i - 1},
    \quad \parbox{2.8cm}{\textit{for} {} $0 < s < 1 $,}
$$
    where $\widehat {F}(t;s) = {{F(t;qs)} \mathord{\left/{\vphantom {{F(t;qs)} q}}
    \right. \kern-\nulldelimiterspace} q}$ is the GF of a sub-critical MBP and
    ${\widehat {F}(t;s)} < 1$, so it converges to one as $i \to \infty $.
    Hence ${\mathcal Q}_{ij} (t) \downarrow 0$ as $i \to \infty $. Last fact
    implies that ${\mathcal Q}_{ij} (t)$ is the Feller function.  Therefore
    $\mathbb{Q}(t)$ has a stable q-matrix with components $q_{ij} = {\mathcal Q}'_{ij}
    (\varepsilon  \downarrow 0)$; see Anderson (1991, p.43).

    Next, since all states are stable then $\mathbb{Q}(t)
    = \left\{ {{\mathcal Q}_{ij} (t)} \right\}$ is  differentiable and has
    a finite and a continuous derivative with respect to $t\in{\mathcal T}$;
    see Anderson (1991, p.10).  Let's compute this derivative.
    It follows from (3.1) that
\begin{eqnarray}
    \Delta {\mathcal Q}_{ij} (t) \nonumber
    & = & {\mathcal Q}_{ij} (t + \varepsilon ) - {\mathcal Q}_{ij} (t) = {{jq^{j - i} }
        \over {i\beta ^t }} \left[ {{{P_{ij} (t + \varepsilon )}
        \over {\beta ^\varepsilon  }} - P_{ij} (t)} \right] \\ \nonumber
    & = & {{jq^{j - i} } \over {i\beta ^t }}\left[ {\Delta P_{ij} (t)
        + P_{ij} (t + \varepsilon ) \ln \beta
        \cdot \varepsilon  + o(\varepsilon )} \right].  \nonumber
\end{eqnarray}
    Hence
$$
    {{\Delta {\mathcal Q}_{ij} (t)} \over \varepsilon } = {{jq^{j - i} }
    \over {i\beta ^t }}\left[ {{{\Delta P_{ij} (t)}
    \over \varepsilon } + P_{ij} (t + \varepsilon )
    \cdot \ln \beta  + o(1)} \right].
$$
    Taking limit as $\varepsilon  \downarrow 0$ here yields
$$
    {\mathcal Q}'_{ij} (t) = {{jq^{j - i} } \over {i\beta ^t }}
    \left[ {P'_{ij} (t) - P_{ij} (t)\ln \beta } \right],
    \quad \parbox{3cm}{\textit{for  all} {} $ i, j \in {\mathcal E} $.}    \eqno (4.8)
$$
    Being that $q_{ij} = {\mathcal Q}'_{ij} (\varepsilon \downarrow 0)$,
    we should compute $P'_{ij} (\varepsilon  \downarrow 0)$.
    It follows from the general theory of MBP that $P'_{ij}
    (\varepsilon  \downarrow 0) = \lim _{\varepsilon  \downarrow 0}
    {{P_{ij} (\varepsilon )} \mathord{\left/ {\vphantom {{P_{ij}
    (\varepsilon )} \varepsilon }} \right. \kern-\nulldelimiterspace}
    \varepsilon } = ia_{j - i + 1}$. Therefore we get that
$$
    q_{ij}  = {{jq^{j - i} } \over i}\left[ {ia_{j - i + 1}
    - \delta _{ij} \ln \beta } \right].
$$
    Using the expression for densities $\left\{ {\lambda _j } \right\}$
    said in (4.2) from the last formula we obtain (4.6). Moreover we see
    that $q_{ij}  \ge 0$ when $i \ne j$ for all $i,j \in {\mathcal E}$ and
    being that both $\lambda _1 $ and $\ln \beta $ are negative yield
    that $q_i : = q_{ii}  < 0$.

    Lastly owing to Markovian nature of $W(t)$ it follows
    from theory of continuous-time Markov chain that
    equation (4.7) holds. In particular, at $\tau = 0$
$$
    {\mathcal Q}'_{ij} (t) = \sum\limits_{k \in {\mathcal E}} {q_{ik} {\mathcal Q}_{kj} (t)}.
$$
    The proof is completed.

    Let $\mathcal{G}_i (s)$ be the GF of q-matrix
    $\left\{ q_{ij} \right\}$  that is
$$
    \mathcal{G}_i (s): = \sum\limits_{j \in {\mathcal E}} {q_{ij} s^j }
    = \sum\limits_{j \in {\mathcal E}} {{\mathcal Q}'_{ij} (\varepsilon  \downarrow 0)s^j }.
$$
    Using expressions (4.6) it follows that
\begin{eqnarray}
    \mathcal{G}_i (s) \nonumber
    & = & \left[ {i\lambda_1  + (i - 1)\ln \beta } \right]s^i
    + \sum\limits_{j \in {\mathcal E}\backslash \{ i\} }
    {{{j\lambda_{j - i + 1} } \over {j - i + 1}}s^j } \\ \nonumber
    & = & (i - 1)s^{i - 1} \left[ {s\ln \beta
    + \sum\limits_{j \in {\mathcal E}} {{{\lambda_j }
    \over j}s^j } } \right] + s^{i - 1} g(s),  \nonumber
\end{eqnarray}
    where $g(s)$ is defined in (4.3). On the other hand it is easy to see that
$$
    \sum\limits_{j \in {\mathcal E}} {{{\, \lambda_j }
    \over j}s^j }  = \int_0^s {{{\, g(u)} \over u}du}.
$$
    Thence we have that
$$
    \mathcal{G}_i (s) = {{\, (i - 1)m(s) + g(s)} \over s}s^i,
$$
    where
$$
    m(s): = s\ln \beta  + \int_0^s {{{\, g(x)} \over x}dx}.
$$

    Now more general, consider
$$
    \mathcal{G}_i (t;s) = \sum\limits_{j \in {\mathcal E}} {{\mathcal Q}'_{ij} (t)s^j }
    = {{\partial G_i (t;s)} \over {\partial t}}\, \raise 1.2pt\hbox{.}
$$
    (The differentiable property of GF $G(t;s)$ will be established in
    the Theorem 14 below). After standard calculations we make sure that
    the GF version of (4.8) is the following identity:
$$
    \mathcal{G}_i (t;s) = {{(i - 1)m\left( {\widehat {F}(t;s)} \right)
    + g\left( {\widehat {F}(t;s)} \right)} \over {\widehat {F}(t;s)}}G_i (t;s),    \eqno (4.9)
$$
    for all $t \in {\mathcal T}$, where $\widehat {F}(t;s) = {{F(t;qs)}
    \mathord{\left/{\vphantom {{F(t;qs)} q}} \right. \kern-\nulldelimiterspace} q}$.

\textbf{Theorem~14.}
    \textit{The GF $G(t;s)$ is differentiable function with respect to
    $t \in {\mathcal T}$ uniformly for $0 \le s < 1$. The transition
    function $\left\{ {{\mathcal Q}_{ij} (t)} \right\}$ is unique
    solution of the backward Kolmogorov system (4.7),
    which is unique GF solution of equation}
$$
    {{\partial G(t;s)} \over {\partial t}}
    = h\left( {\widehat F(t;s)} \right)G(t;s),      \eqno (4.10)
$$
    \textit{with condition $G(0;s) = s$, where
    $h(s) = {{g(s)} \mathord{\left/ {\vphantom {{g(s)} s}}
    \right. \kern-\nulldelimiterspace} s}$.}

\textbf{\textit{Proof}.}
    As in proof of the Theorem 13 we have
$$
    {\mathcal Q}_{ij} (t) - {\mathcal Q}_{ij} (t + \varepsilon ) \le
    {\mathcal Q}_{ij} (t) \cdot \left[ {1 - {\mathcal Q}_{ii} (\varepsilon )} \right]
$$
    for arbitrary $\varepsilon  > 0$. Hence for the difference
$$
    \Delta _\varepsilon  G(t;s) = G(t + \varepsilon ;s) - G(t;s)
$$
    we obtain that
\begin{eqnarray}
    \left|\Delta _\varepsilon  G(t;s)\right| \nonumber
    & \le & \sum\limits_{j \in {\mathcal E}} {\left| {{\mathcal Q}_{1j} (t)
        - {\mathcal Q}_{1j} (t + \varepsilon )} \right|s^j } \\ \nonumber
    & \le & \left[ {1 - {\mathcal Q}_{11} (\varepsilon )} \right] \cdot
        \sum\limits_{j \in {\mathcal E}} {{{\mathcal Q}_{1j}(t)}s^j} = 2G(t;s)
    \cdot \left[ {1 - {\mathcal Q}_{11} (\varepsilon )} \right]. \nonumber
\end{eqnarray}
    Since ${\mathcal Q}_{ij} (t)$ is standard, it follows from
    last inequality that $\Delta _\varepsilon  G(t;s) \to 0$ as
    $\varepsilon  \downarrow 0$. So $G(t;s)$ is continuous function with
    respect to $t \in {\mathcal T}$ uniformly for $0 \le s < 1$. It can
    easily be seen that a GF version of the relation (4.2) is
$$
    G(\varepsilon ;s) = s + g(s) \cdot \varepsilon  + o(\varepsilon )
    \quad \parbox{2cm}{\textit{as} {} $\varepsilon  \downarrow 0$,}   \eqno (4.11)
$$
    for $0 \le s < 1$. By the way according to formulas (1.4) and (4.5)
    one can see that GF $G(t;s)$ satisfies the following functional equation:
$$
    G(t + \tau ;s) = \,{{G\left( {t;\widehat F(\tau ;s)} \right)}
    \over {\,G\left( {0;\widehat F(\tau ;s)} \right)\,}}G(\tau ;s).    \eqno (4.12)
$$
    We use expressions (4.11) and (4.12) to
    $\Delta _\varepsilon  G(t;s)$ and hereupon we get
$$
    \Delta _\varepsilon G(t;s) = {{g\left({\widehat F(t;s)} \right)}
    \over {\widehat F(t;s)}}  \cdot \varepsilon + o(\varepsilon )
    \quad \parbox{2cm}{\textit{as} {} $\varepsilon  \downarrow 0$,}
$$
    for any $t \in {\mathcal T}$ and all $0 \le s < 1$,
    which implies that $G(t;s)$ is differentiable.

    Equation (4.10) follows from formula (4.9) at $i = 1$ and
    the boundary condition $G(0;s) = s$ follows from (4.11).
    The uniqueness of solution of (4.10) follows from the
    classical differential equations theory.

    The theorem is proved.

    The following assertion is direct consequence from Theorem 14.

\textbf{Corollary.}
    \textit{The differential equation (4.10) is equivalent to the following one}
$$
    \int_0^t {h\left( {\widehat F(\tau ;s)}
    \right)d\tau }  = \ln {{G(t;s)} \over s}   \eqno (4.13)
$$
    \textit{with boundary condition $G(0;s) = s$, where $h(s) = {{g(s)}
    \mathord{\left/ {\vphantom {{g(s)} s}}\right. \kern-\nulldelimiterspace} s}$}.

\subsection{Classification and Ergodic behavior}

    As it has been noticed above, that the parameter $a = f'(s \uparrow 1)$
    plays a regulating role for MBP and is subdivided three types of process
    depending on sign of $a$. Note that evolution of MQP is regulated in
    essence by positive parameter $\beta  = \exp \{ f'(q)\} $. Thus are
    subdivided two types of process depending on values of this parameter.
    From equalities (4.5) and (4.13) we write
$$
    G_i (t;s) = s\left[ {\widehat F(t;s)} \right]^{i - 1} \exp
    \left\{ {\int_0^t {h\left( {\widehat F(\tau ;s)} \right)d\tau } } \right\}.   \eqno (4.14)
$$

    If $\alpha : = g'(1)$ is finite, then it follows from (1.14) that
$$
    \mathbb{E}_i W(t) = \left({i - 1}\right)\beta ^t + \mathbb{E} W(t)
$$
    and
\begin{equation}
    \mathbb{E} W(t)= \left\{\begin{array}{l} 1 + \gamma \left( {1 - \beta ^t }
    \right) \hfill, \quad \parbox{2.5cm}{\textit{when} {} $\beta  < 1 $,}    \nonumber\\
    \\
    \alpha t + 1 \hfill,  \quad \parbox{2.5cm}{\textit{when} {} $\beta  = 1 $,}   \nonumber\\
    \end{array} \right.    \eqno (4.15)
\end{equation}
    where $\gamma  = {\alpha  \mathord{\left/ {\vphantom {\alpha
    {\left| {\ln \beta } \right|}}} \right. \kern-\nulldelimiterspace}
    {\left| {\ln \beta } \right|}}$. Moreover we obtain the variance structure
\begin{equation}
    \textsf{Var}_i W(t) = \left\{\begin{array}{l} \left[ {\gamma
    + \left( {i - 1} \right)\left( {1 + \gamma } \right)\beta ^t }
    \right]\left( {1 - \beta ^t } \right) \hfill,
    \quad \parbox{2.5cm}{\textit{when} {} $\beta  < 1 $,}   \nonumber \\
    \\
    \alpha it \hfill, \quad \parbox{2.5cm}{\textit{when} {} $\beta =1 $,} \nonumber\\
    \end{array} \right.
\end{equation}
     where  $\textsf{Var}_i W(t) = \textsf{Var}
     \left[ {W(t)\left| {W(0) = i} \right.} \right]$.

    The formula (4.15) implies that when $\beta  = 1$
$$
    \mathbb{E}_i W(t) \sim \alpha t
    \quad \parbox{2.2cm}{\textit{as} {} $ t \to \infty  $,}
$$
    and if $0 < \beta  < 1$ then
$$
    \mathbb{E}_i W(t) \longrightarrow 1 + \gamma
    \quad \parbox{2.2cm}{\textit{as} {} $ t \to \infty $.}
$$
    Thereby we classify the MQP as \textit{restrictive}
    if $\beta  < 1$ and \textit{explosive} if $\beta  = 1$.

    Further in restrictive case we keep on the condition (2.12).
    As it was noted in Section 2, this condition is equivalently to
    $\sum\nolimits_{j \in {\mathcal S}} {a_j q^{j - 1} j\ln j}  < \infty$.
    Being that $\lambda _j  = jq^{j - 1} a_j $, for feasibility of
    (2.12) it is necessary and sufficient that
$$
    \sum\limits_{j \in {\mathcal E}} {\lambda _j \ln j}  < \infty.   \eqno (4.16)
$$
    In explosive case we everywhere suppose that the condition (1.7) is satisfied.

\textbf{Theorem~15.}
    \textit{The MQP is}
\begin{enumerate}
    \item [\textbf{(i)}]   \textit{positive if it is
                            restrictive and condition (4.16) is satisfied;}

    \item [\textbf{(ii)}]  \textit{null if it is explosive.}
\end{enumerate}

\textbf{\textit{Proof}.}
    To prove assertion (i) from (4.13) we get
$$
    \ln {\mathcal Q}_{11} (t) = \int_0^t {h\left( {\widehat F(\tau ;0)} \right)d\tau }
    = \int_0^{\widehat F(t;0)} {{{h(x)} \over {\widehat f(x)}}dx}
    \longrightarrow \int_0^1 {{{h(x)} \over {\widehat f(x)}}dx},
$$
    since $\widehat F(t;0) \uparrow 1$ as $t \to \infty $, where
    $\widehat f(s) = {{f(qs)} \mathord{\left/ {\vphantom {{f(qs)} q}}
    \right. \kern-\nulldelimiterspace} q}$. With reference
    to Yang (1972) we make sure that the condition (4.16) is
    sufficient for a converging the integral in right-hand side.
    Hence $\lim _{t \to \infty } {\mathcal Q}_{11} (t) > 0$. For part (ii)
    we recall that $q = 1$ and ${{h(s)} \mathord{\left/ {\vphantom {{h(s)} s}}
    \right. \kern-\nulldelimiterspace} s} = f'(s)$ if $\beta  = 1$. Similarly
$$
    \ln {\mathcal Q}_{11}(t)= \int_0^t {h\left( {F(\tau ;0)}\right)d\tau }
    = \int_0^{F(t;0)} {{{h(x)} \over {f(x)}}dx}  \longrightarrow
    \int_0^1 {{{f'(x)} \over {f(x)}}dx}  =  - \infty,
$$
    as $t \to \infty $, so that $\lim _{t \to \infty } {\mathcal Q}_{11} (t) = 0$.

    The theorem is proved.

    The next two assertions are direct consequences of Lemma 2.

\textbf{Theorem~16.}
    \textit{Let MQP be restrictive. If condition (4.16) is
    satisfied, then }
$$
    G_i (t;s) = {\mathcal U}(s)\left( {1 + o(1)} \right)
    \quad \parbox{2.2cm}{\textit{as} {} $ t \to \infty $,}
$$
    \textit{for all $0 \le s < 1$, where the limiting
    GF ${\mathcal U}(s) = \sum\nolimits_{j \in {\mathcal E}} {u_j s^j } $ has a form}
$$
    {\mathcal U}(s) = s{{|\ln \beta |} \over {f(qs)}}{\mathcal A}(qs).  \eqno (4.17)
$$
    \textit{The numbers $\left\{ {u_j } \right\}$ represent an invariant distribution for MQP.}

\textbf{\textit{Proof}.}
    The convergence of $G_i (t;s)$ to ${\mathcal U}(s)$ follows from assertion (2.13) and
    formula (4.5) because $\widehat F(t;s) \uparrow 1$ as $t \to \infty $,
    where ${\mathcal U}(s)$ in the form of (4.17). Taking limit in (4.12)
    implies a Schroeder type invariance equation
$$
    {\mathcal U}(s) = {{G(\tau ;s)} \over {\,\widehat F(\tau ;s)\,}}
    {\mathcal U}\left( {\widehat F(\tau ;s)} \right)
$$
    and hence $u_j  = \sum\nolimits_{i \in {\mathcal E}} {u_i {\mathcal Q}_{ij} (\tau )} $
    for any $\tau  \in {\mathcal T}$. Therefore $\left\{ {u_j } \right\}$ is an
    invariant measure. Let now condition (4.16) be satisfied.
    Then according to properties of the function ${\mathcal A}(s)$ (see (2.10))
$$
    \sum\limits_{j \in {\mathcal E}} {u_j }  = \mathop {\lim }\limits_{s \uparrow 1}
    {\mathcal U}(s) = \mathop {\lim }\limits_{s \uparrow 1} {{{\mathcal A}(qs)}
    \over {q\left( {1 - s} \right)}} = 1.
$$
    The theorem is proved.

\textbf{Theorem~17.}
    \textit{If MQP is explosive, then for all $0 \le s < 1$}
$$
    {{(\nu t)^{1 + {1 \mathord{\left/ {\vphantom {1 \nu }} \right.
    \kern-\nulldelimiterspace} \nu}}}\over {{\mathcal N}(t)}}\cdot G_i (t;s)
    \longrightarrow \pi (s)
    \quad \parbox{2.2cm}{\textit{as} {} $ t \to \infty $,}       \eqno (4.18)
$$
    \textit{where ${\mathcal N}(t)$ satisfies the property (2.3).
    The limit GF $\pi (s) = \sum\nolimits_{j \in {\mathcal E}} {\pi _j s^j } $
    determines an invariant measure $\left\{ {\pi _j } \right\}$ and}
$$
    \sum\limits_{j = 1}^n {\pi _j }  = {1 \over
    {\Gamma (2 + \nu )}}n^{1 + \nu } \mathfrak{L}_{\pi } (n),     \eqno (4.19)
$$
    \textit{where $\Gamma ( * )$ is Euler's Gamma function and
    $\mathfrak{L}_{\pi }(n) \cdot \mathfrak{L}(n) \to 1$ as $n \to \infty$.}

\textbf{\textit{Proof}.}
    From second part of Lemma 2 and (4.5) we will write out
$$
    G_i (t;s) \sim s{{R^{1 + \nu } (t;s)} \over {f(s)}} \cdot
    \mathfrak{L}\left( {{1 \over {R(t;s)}}} \right)
    \quad \parbox{2.2cm}{\textit{as} {} $ t \to \infty $.}     \eqno (4.20)
$$
    It follows from (1.7) and (2.2) that
$$
    {{R^{1 + \nu } (t;s)} \over {f(s)}} = {{{\mathcal N}^{1 + \nu } (t)} \over
    {(\nu t)^{1 + {1 \mathord{\left/ {\vphantom {1 \nu }} \right.
    \kern-\nulldelimiterspace} \nu }} }} \cdot {{\mathfrak{L}^{- 1}
    \left( {{1 \mathord{\left/ {\vphantom {1 {(1 - s)}}} \right.
    \kern-\nulldelimiterspace} {(1 - s)}}} \right)} \over {(1 - s)^{1 + \nu } }}
    \cdot {1 \over {\left( {1 + {{\mathcal{M}(s)} \mathord{\left/{\vphantom
    {{\mathcal{M}(s)} t}} \right. \kern-\nulldelimiterspace} t}} \right)^{1
    + {1 \mathord{\left/ {\vphantom {1 \nu }} \right.
    \kern-\nulldelimiterspace} \nu }} }},               \eqno (4.21)
$$
    where the function $\mathcal{M}(s)$ is defined in (2.4). It is cogently
    that $\mathcal{M}(0) = 0$ and ${{R(t;s)} \mathord{\left/ {\vphantom
    {{R(t;s)} {R(t;0)}}} \right. \kern-\nulldelimiterspace} {R(t;0)}} \to 1$
    as $t \to \infty $ uniformly for $0 \le s < 1$. Hence according to (2.3),
    ${\mathcal N}^{\nu }(t) \cdot \mathfrak{L}\left({{1 \mathord{\left/
    {\vphantom {1 {R(t;s)}}} \right. \kern-\nulldelimiterspace}
    {R(t;s)}}} \right) \to 1$. Then from (4.20) and (4.21) appears
$$
    G_i (t;s) \sim {{{\mathcal N}(t)} \over {(\nu t)^{1 + {1 \mathord{\left/ {\vphantom
    {1 \nu }} \right. \kern-\nulldelimiterspace} \nu }} }} \cdot \pi (s) \cdot
    {1 \over {\left( {1 + {{\mathcal{M}(s)} \mathord{\left/ {\vphantom
    {{\mathcal{M}(s)} t}} \right. \kern-\nulldelimiterspace} t}} \right)^{1 +
    {1 \mathord{\left/ {\vphantom {1 \nu }} \right.
    \kern-\nulldelimiterspace} \nu }} }},            \eqno (4.22)
$$
    as $t \to \infty$, where
$$
    \pi (s) = {s \over {(1 - s)^{1 + \nu } }}\mathfrak{L}_{\pi }
    \left( {{1 \over {1 - s}}} \right).         \eqno (4.23)
$$
    The expansion (4.18) follows from the relation (4.22). The invariant equation
    $\pi _j  = \sum\nolimits_{i \in {\mathcal E}} {\pi _i {\mathcal Q}_{ij} (t)} $ comes
    out from the functional equation (4.12). At last, according to the Hardy-Littlewood
    Tauberian theorem each of relations (4.19) and (4.23) entails another.

    The theorem is proved.

    It undoubtedly that $\lim _{s \downarrow 0} \left[ {{{G_i (t;s)}
    \mathord{\left/ {\vphantom {{G_i (t;s)} s}} \right. \kern-\nulldelimiterspace} s}}
    \right] = {\mathcal Q}_{i1} (t)$. Then from Theorems 16 and 17 we get to
    the following local limit theorems.

\textbf{Theorem~18.}
    \textit{If MQP is restrictive and condition (4.16) is satisfied, then}
$$
    {\mathcal Q}_{i1} (t) \longrightarrow{{\left| {\ln \beta } \right|}
    \over {a_0 }}{\mathcal A}(0)
    \quad \parbox{2.2cm}{\textit{as} {} $ t \to \infty $.}
$$

\textbf{Theorem~19.}
    \textit{If MQP is explosive, then}
$$
    (\nu t)^{1 + {1 \mathord{\left/ {\vphantom {1 \nu }} \right.
    \kern-\nulldelimiterspace} \nu }}\cdot Q_{i1} (t)
    \sim {{{\mathcal N} (t)} \over {a_0 }}
    \quad \parbox{2cm}{\textit{as} {} $t \to \infty $,}    \eqno (2.16)
$$
    \textit{where the function ${\mathcal N}(t)$ satisfies the property (2.3).}

    Further we observe limit properties of $\left\{ {{{{\mathcal Q}_{ij} (t)} \mathord
    {\left/{\vphantom {{{\mathcal Q}_{ij} (t)} {{\mathcal Q}_{11} (t)}}} \right.
    \kern-\nulldelimiterspace} {{\mathcal Q}_{11} (t)}}} \right\}$. Consider the GF
$$
    {\mathcal W}_i (t;s) = \sum\limits_{j \in {\mathcal E}} {{{{\mathcal Q}_{ij} (t)} \over
    {{\mathcal Q}_{11} (t)}}s^j }  = {1 \over {{\mathcal Q}_{11} (t)}}G_i (t;s)
    = \left[ {\widehat F(t;s)} \right]^{i - 1} {\mathcal W}(t;s),       \eqno (4.24)
$$
    where
$$
    {\mathcal W}(t;s) = \sum\limits_{j \in {\mathcal E}}
    {{{{\mathcal Q}_{1j} (t)} \over {{\mathcal Q}_{11} (t)}}s^j }.
$$
    For general MQP the following ratio limit property holds.

\textbf{Theorem~20.}
    \textit{The limits}
$$
    \mathop {\lim }\limits_{t \to \infty } {{{\mathcal Q}_{ij} (t)}
    \over {{\mathcal Q}_{11} (t)}} =: \omega _j        \eqno (4.25)
$$
    \textit{exist for all $i,j \in {\mathcal E}$. The set
    $\left\{ {\omega _j } \right\}$ is an invariant
    measure and the GF}
$$
    {\mathcal U}(s) := \sum\limits_{j \in {\mathcal E}} {\omega _j s^j }
    = s\,\exp \left\{ {\int_0^s {{{|h(x)|} \over
    {\widehat f(x)}}dx} } \right\},        \eqno (4.26)
$$
    \textit{converges for $0 \le s < 1$, where $h(s) = {{g(s)}
    \mathord{\left/ {\vphantom {{g(s)} s}} \right. \kern-\nulldelimiterspace} s}$
    and  $\widehat f(s) = {{f(qs)} \mathord{\left/ {\vphantom {{f(qs)} q}} \right.
    \kern-\nulldelimiterspace} q}$.}

\textbf{\textit{Proof}.}
    It follows from (4.24) that it suffice to consider the case $i = 1$
    because $\widehat F(t;s) \uparrow 1$ as $t \to \infty $ uniformly
    for all $0 \le s \leq r < 1$. So write
$$
    {\mathcal U}(t;s) = s \exp \left\{ {\int_0^t {\left[ {h\left( {\widehat F(u;s)} \right)
    - h\left( {\widehat F(u;0)} \right)} \right]du} } \right\}.
$$
    One can choose $\tau  \in {\mathcal T}$ for any $0 \le s < 1$ so that
    $s = \widehat F(\tau ;0)$. Then considering functional equation (1.4),
    we get  $\widehat F( {t;s} ) = \widehat F(t + \tau ;0)$ and hence
\begin{eqnarray}
    {\mathcal U}(t;s) \nonumber
    & = & s \exp \left\{ {\int_\tau ^{t + \tau } {h\left( {\widehat F(u;0)} \right)du}
    - \int_0^t {h\left( {\widehat F(u;0)} \right)du} } \right\} \\   \nonumber
    \\  \nonumber
    & = & s \exp \left\{ {\int_0^\tau  {\left[ {h\left( {\widehat F(t;\widehat F(u;0))} \right)
    - h\left( {\widehat F(u;0)} \right)} \right]du} } \right\} \\   \nonumber
    \\   \nonumber
    & = & s \exp \left\{ {\int_0^s {{{h\left( {\widehat F(t;x)} \right)
    - h(x)} \over {\widehat f(x)}}dx} } \right\},  \nonumber
\end{eqnarray}
    where we have used the equation (1.5) and $\widehat f(s) =
    {{f(qs)} \mathord{\left/ {\vphantom {{f(qs)} q}} \right. \kern-\nulldelimiterspace} q}$.
    To get to (4.26) it suffice to take limit as $t \to \infty $ being that
    $\widehat F(t;s) \to 1$ and $h(1) = 0$. Assertion (4.25) follows from (4.26) owing
    to the continuity theorem for GF. It easily to be convinced
    ${\mathcal W}(s) < \infty $ for all $0 \le s < 1$.

    Now we observe that the set $\left\{ {\omega _j } \right\}$ to be the invariant
    measure for MQP. First using the Kolmogorov-Chapman equation we obtain that
$$
    {{{\mathcal Q}_{ij} (t + \tau )} \over {{\mathcal Q}_{11} (t + \tau )}} \cdot
    {{{\mathcal Q}_{11} (t + \tau )} \over {{\mathcal Q}_{11} (t)}} = \sum\limits_{k \in {\mathcal E}}
    {{{{\mathcal Q}_{ik} (t)} \over {{\mathcal Q}_{11} (t)}}{\mathcal Q}_{kj} (\tau )}.
$$
    Setting $s = 0$ in (4.12) we can see that ${{{\mathcal Q}_{11} (t + \tau )}
    \mathord{\left/ {\vphantom {{{\mathcal Q}_{11} (t + \tau )} {{\mathcal Q}_{11} (t)}}} \right.
    \kern-\nulldelimiterspace} {{\mathcal Q}_{11} (t)}} \to 1$ as $t \to \infty $. Hence we
    get the following invariant equation for $\left\{ {\omega _j } \right\}$
$$
    \omega _j  = \sum\limits_{k \in {\mathcal E}} {\omega _k {\mathcal Q}_{kj} (t)}
    \quad \parbox{3cm}{\textit{for \,any} {} $ t \in {\mathcal T} $.}        \eqno (4.27)
$$
    The GF version of (4.27) is
$$
    {\mathcal W}\left( {\widehat F(t;s)} \right)
    = {{\widehat F(t;s)} \over {\,G(t;s)\,}}{\mathcal W}(s),
$$
    for $0 \le s < 1$, the functional equation of generalized Schroeder form.

    The theorem is proved.

    We complete the paper with stating of the following limit theorem.

\textbf{Theorem~21.}
    \textit{Let MQP is explosive and the function ${\mathcal N}(t)$ satisfies
    the property (2.3). Then for any $x > 0$}
$$
    \mathbb{P}_i \left\{ {{{{\mathcal N}(t)} \over {(\nu n)^{{1 \mathord{\left/ {\vphantom
    {1 \nu }} \right. \kern-\nulldelimiterspace} \nu }} }}W(t) < x} \right\}
    \longrightarrow G(x)
    \quad \parbox{2.2cm}{\textit{as} {} $ t \to \infty $,}
$$
    \textit{where the Laplace transform}
$$
    \int_{\mathbb{R}_+} {e^{ - \theta x} dG(x)}  = {1 \over {\left( {1
    + \theta ^\nu  } \right)^{1 + {1 \mathord{\left/ {\vphantom {1 \nu }}
    \right. \kern-\nulldelimiterspace} \nu }} }}.
$$

\textbf{\textit{Proof}.}
    Consider the Laplace transform
$$
    \phi (t;\theta ): = \mathbb{E}e^{ - \theta q(t)W(t)}
    = G_i \left( {t;\theta (t)} \right),         \eqno (4.28)
$$
    where $\theta (t) = \exp \{  - \theta q(t)\} $ and
$$
    q(t): = R\left( {t;0} \right) = {{{\mathcal N}(t)} \over {(\nu t)^{{1 \mathord
    {\left/ {\vphantom {1 \nu }} \right. \kern-\nulldelimiterspace} \nu }} }}.
$$
    It was shown in the proof of Theorem 17 that
$$
    G_i (t;s) \sim {{q(t)} \over {\nu t}} \cdot \pi (s) \cdot {1 \over
    {\left( {1 + {{\mathcal{M}(s)} \mathord{\left/ {\vphantom {{M(s)} t}}
    \right. \kern-\nulldelimiterspace} t}} \right)^{1 + {1 \mathord
    {\left/ {\vphantom {1 \nu }} \right. \kern-\nulldelimiterspace} \nu }}}}
    \quad \parbox{2.2cm}{\textit{as} {} $ t \to \infty $,}       \eqno (4.29)
$$
    where
$$
    \pi (s) = {s \over {(1 - s)^{1 + \nu } \mathfrak{L}\left( {{1
    \mathord{\left/ {\vphantom {1 {(1 - s)}}} \right.
    \kern-\nulldelimiterspace} {(1 - s)}}} \right)}}    \eqno (4.30)
$$
    and it follows from (3.5) and (3.8) that
$$
    \mathcal{M}(s) \sim {1 \over {\nu (1 - s)^\nu  \mathfrak{L}\left( {{1
    \mathord{\left/ {\vphantom {1 {(1 - s)}}} \right.
    \kern-\nulldelimiterspace} {(1 - s)}}} \right)}}
    \quad \parbox{2.2cm}{\textit{as} {} $ s \uparrow 1 $.}       \eqno (4.31)
$$
    Put $s = \theta (t)$ in (4.29). It is clear $1 - \theta (t) \sim \theta q(t)$ and
    ${{R(t;s)} \mathord{\left/ {\vphantom {{R(t;s)} {q(t)}}} \right.
    \kern-\nulldelimiterspace} {q(t)}} \to 1$ as $t \to \infty $. So by the
    property of SV functions $\mathfrak{L}\left( {{1 \mathord{\left/
    {\vphantom {1 {R\left( {t;\theta (t)} \right)}}} \right. \kern-\nulldelimiterspace}
    {R\left( {t;\theta (t)} \right)}}} \right) \sim \mathfrak{L}\left( {{1
    \mathord{\left/{\vphantom {1 {\left( {1 - \theta (t)} \right)}}} \right.
    \kern-\nulldelimiterspace} {\left( {1 - \theta (t)} \right)}}} \right)$.
    On account of all these, from (4.28)--(4.31) it is a matter of standard
    computation to verify that
$$
    \phi (t;\theta ) \longrightarrow  {1 \over {\left( {1 + \theta ^\nu  }
    \right)^{1 + {1 \mathord{\left/ {\vphantom {1 \nu }}  \right.
    \kern-\nulldelimiterspace} \nu }} }}
    \quad \parbox{2.2cm}{\textit{as} {} $ t \to \infty $.}
$$

    The theorem proof is completed.

    The Theorem 21 generalizes for all $0 < \nu  \le 1$ the known Harris
    theorem established under finite variance condition for the process
    with discrete time; see Athreya and Ney (1972, p.59). Really, in
    specific case $\nu  = 1$, the Laplace transform specified in the
    theorem becomes  $\left( {1 + \theta } \right)^{ - 2} $ that
    fits to the first order Erlang's law
$$
    1 - e^{ - x}  - xe^{ - x}.
$$

\section{Conclusion}

    We devote the paper to research of the asymptote of trajectory and limit structure
    of MBP $\left\{ {Z(t), t\in {\mathcal T}} \right\}$. All our reasoning and results are
    based on the assertion of the Lemma 1 and Lemma 2. Local limit theorems proved in
    Section 2 and Section 3 improve same results from the  paper of Imomov (2014a)
    excepting a finiteness of the second moment $f''(s\uparrow 1)$.

    In the Section 4 we observe the structural and asymptotical properties of
    Markov Q-process (MQP). This process represents a stochastic population
    process with the trajectory never extinct. We see that MQP is subdivided
    two types of process depending on values of parameter $\beta =\exp \{f'(q)\}$.
    Classification and ergodic properties of MQP are studied in this Section.

    In our next researches we will improve the results for critical MBP provided
    that  $\mathfrak{L}(\cdot)$ is the normalized slowly varying function with remainder so that
$$
    {{\mathfrak{L}\left( {\lambda x} \right)} \over {\mathfrak{L}(x)}}
    = 1 + {o}\left( {{{\mathfrak{L}\left(x\right)} \over {x^\nu }}} \right)
    \quad \parbox{2.2cm}{\textit{as} {} $x  \rightarrow \infty$}
$$
    for each $\lambda \in \mathbb{R}_+ $; see Bingham \textit{et al.} (1987, p.185).
    We already have some advancement at this assumption for a discrete case only.
    So if GF $F(s)$ of offspring law of discrete time branching process
    $\left\{ {Z(n), n\in  \mathbb{N}} \right\}$ has a representation
$$
    F(s) = s + (1 - s)^{1 + \nu } \mathfrak{L}
    \left( {{1 \over {1 - s}}} \right),
$$
    then
$$
    {1 \over {\Lambda \bigl( {1-F(n;s)} \bigr)}} - {1 \over {\Lambda \left
    ({1-s} \right)}} \sim \nu n + {{1+\nu }\over 2}\cdot
    \ln \left( 1+ \nu n \Lambda (1-s)\right)             \eqno (5.1)
$$
    as $n \to \infty$, where $F(n;s)=\mathbb{E}{s^{Z(n)}}$ and
    $\Lambda (y) = y^\nu {\mathfrak L}\left( {{1 \mathord{\left/
    {\vphantom {1 y}} \right. \kern-\nulldelimiterspace} y}} \right)$.

    We are sure the statement (5.1) is fair and for MBP. Hence we can state convergence
    rates in limit theorems in the case $a=0$ for MBP and in the case $\beta =1$ for MQP.

\bigskip

\begin{center}
\textbf{ \Large References}
\end{center}

{\leftskip=1em\parindent=-1em

Anderson W. (1991). \textit{Continuous-Time Markov Chains: An Applications-Oriented Approach}. Springer, New York.

Asmussen S., Hering H. (1983). \textit{Branching Processes}. Boston.

Athreya K.B. and Ney P.E. (1972). \textit{Branching processes}. Springer, New York.

Bingham N.H., Goldie C.M., Teugels J.L. (1987). \textit{Regular Variation}. Cambridge.

Bharucha-Reid A.T. (1960). \textit{Elements of the theory of Markov processes and their applications}. New York.

D'Ancona U. (1954). \textit{The struggle for existence}. Bibliotheca Biotheoretica, series D, Leiden, 6.

Feller W. (1939). Die Grundlagen der Volterraschen Theorie des Kampfes ums Dasien in wahrscheinlichkeitstheoretischer Behandlung. \textit{Acta Biometrica}, 5:11--40.

Heatcote C.R., Seneta E. and Vere-Jones. (1967). A refinement of two theorems in the theory of branching process. \textit{Theory of Probability and its Applications}, 12(2):341--346.

Imomov A.A. (2014a). Limit properties of transition functions of continuous-time Markov branching processes. \textit{Intern. Jour. Stochastic Analysis}, 2014:10 pages, DOI: dx/doi.org/10.1155/2014/409345.

Imomov A.A. (2014b). On long-term behavior of continuous-time Markov Branching Processes allowing Immigration. \textit{Journal of Siberian Federal University. Mathematics and Physics}, 7(4):429--440.

Imomov A.A. (2014c). Limit Theorem for the Joint Distribution in the Q-processes. \textit{Journal of Siberian Federal University. Mathematics and Physics}, 7(3):289--296.

Imomov A.A. (2012). On Markov analogue of Q-processes with continuous time. \textit{Theory of Probability and Mathematical Statistics}, 84:57--64.

Imomov A.A. (2002). Some asymptotical behaviors of Galton-Watson branching processes under condition of non-extinctinity of it remote future. \textit{Abstracts of Comm. of VIII Vilnius Conference: Probab. Theory and Math. Statistics}, Vilnius, Lithuania, 118.

Imomov A.A. (2001). On a form of condition of non-extinction of branching processes. \textit{Uzbek Mathematical Journal}, 2:46--51. (Russian)

Karamata J. (1933). Sur un mode de croissance reguliere. Theoremes fondamenteaux. \textit{Bull. Soc. Math. France}, 61:55--62.

Kendall D.G. (1948a). On the generalized "birth-and-death" process. \textit{The Annals of Mathematical Statistics}, 19:1--15.

Kendall D.G. (1948b). On the role of variable generation time in the development of a stochastic birth process. \textit{Biometrica}, 35:316--330.

Kolmogorov A.N and Dmitriev N.A. (1947). Branching stochastic process. \textit{Reports of Academy of Sciences of USSR}, 56:7--10. (Russian)

Lamperti J. and Ney P.E. (1968). Conditioned branching processes and their limiting diffusions. \textit{Theory of Probability and its Applications}, 13:126--137.

Li J., Chen A. and Pakes A.G. (2010). Asymptotic properties of the Markov Branching Process with Immigration. \textit{Journal of Theoretical Probability}, 25(1):122--143.

Nagaev A.V and Badalbaev I.S. (1967). A refinement of certain theorems on branching random process. \textit{Litovskiy Matematicheskiy Sbornik}, 7(1):129--136. (Russian)

Neyman J. (1961). Contributions to biology and problems of medicine. \textit{4-d Berkley Symposium}, 4.

Neyman J. (1956). Contributions to biology and problems of health. \textit{3-d Berkley Symposium}, 4.

Pakes A.G. (2010). Critical Markov branching process limit theorems allowing infinite variance. \textit{Advances in Applied Probability}, 42:460--488.

Pakes A.G. (1999). Revisiting conditional limit theorems for the mortal simple branching process. \textit{Bernoulli}, 5(6):969--998.

Pakes A.G. (1971). Some limit theorems for the total progeny of a branching process. \textit{Advances in Applied Probability}, 3:176--192.

Sevastyanov B.A. (1951). The theory of Branching stochastic process. \textit{Uspehi Mat. Nauk.}, 6(46):47--99. (Russian)

Yang Y.S. (1972). On branching processes allowing immigration. \textit{Journal of Applied Prob.}, 9:24--31.

Yule G.U. (1924). A mathematical theory of evolution based on the conclusions of Dr. J.C.Wills, F.R.S. \textit{Philosophical Transactions of the Royal Society of London}, B213:21--87.

Zolotarev V.M. (1957). More exact statements of several theorems in the theory of branching processes. \textit{Theory of Probability and its Applications}, 2:245--253.

Watson H. and Galton F. (1874). On the probability of the extinction of families. \textit{Journal of Anthropol. Inst. Great Britain and Ireland}, 4:138--144.

}

\end{document}